\documentclass[fleqn,11pt]{article}
\usepackage{graphicx}
\usepackage{amscd,amssymb}
\usepackage{amsmath}

\title{Structured
Vector Bundles Define\\Differential $K$-Theory}
\author{James Simons and Dennis Sullivan}
\date{}

\begin{document}

\maketitle

\begin{abstract}
A equivalence relation, preserving the Chern-Weil form, is defined
between connections on a complex vector bundle.  Bundles equipped
with such an equivalence class are called Structured Bundles, and
their isomorphism classes form an abelian semi-ring.  By applying the
Grothedieck construction one obtains the ring $\hat{K}$, elements of which,
modulo a complex torus of dimension the sum of the odd Betti numbers of
the base, are uniquely determined by the corresponding element of
ordinary $K$ and the Chern-Weil form.  This construction provides a
simple model of differential $K$-theory, c.f. Hopkins-Singer (2005), as
well as a useful codification of vector bundles with connection.
\end{abstract}

\section*{Introduction}
This paper grew out of the effort to construct a simple geometric
model for differential $K$-theory, the fibre
product of usual $K$-theory with closed differential forms, [4],[5],[6].  The model which
finally emerged also fulfilled our long standing wish for a simple and
straightforward codification of complex vector bundles with connection.

Considering pairs of connections whose Chern-Simons
difference form is exact defines an equivalence relation in the space of
all connections on a given bundle.  We call a pair, $\mathcal{V} = (V,\{\nabla\})$,
consisting of a vector bundle together with
a particular such equivalence class, a \textbf{structured bundle}.  As is
true for vector bundles, structured bundles have additive inverses up
to trivial structured bundles: given $\mathcal{V}$ there is a $\mathcal{W}$ such that their direct
sum is equivalent to a bundle with trivial holonomy (Theorem 1.15).

By defining Struct to be the commutative semi-ring of isomorphism
classes of structured bundles, and using the standard Grothedieck
device to turn Struct into a commutative ring, we obtain $\hat{K}$, a
functor from smooth compact manifolds with corners into commutative
rings.  As in ordinary $K$, every element of $\hat{K}$ may be written
as $\mathcal{V} - [n]$, where $\mathcal{V}$ is a structured bundle and $[n]$ is the
trivial structured bundle of dim $n$. $\hat{K}$ achieves the above
desired codification of connections and serves as the sought after 
geometric model of differential $K$-theory.  

Defining four natural transformations into and out of $\hat{K}$ we
develop in the first four sections the diagram with exact diagonals and boundaries,

\vspace{.5cm}
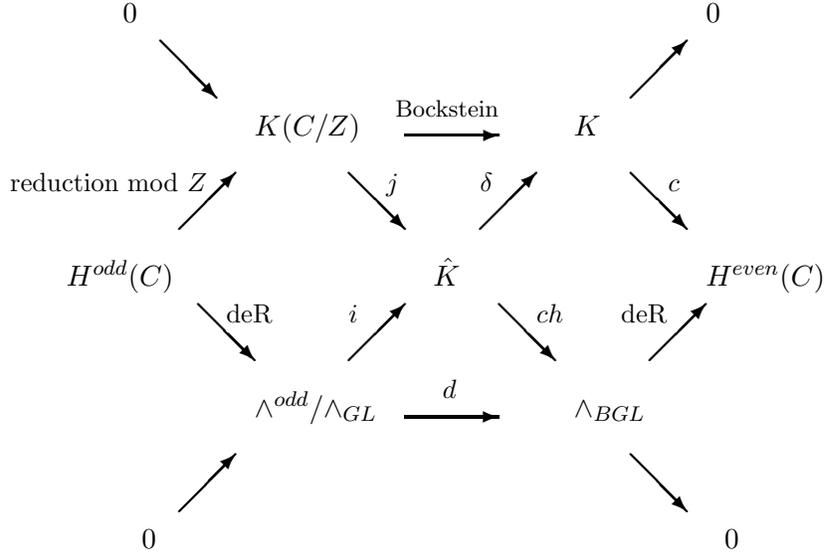
\begin{center}
\setlength{\unitlength}{0.5cm}
\begin{picture}(24,16)\thicklines
\put(5,1){$0$}
\put(20.5,1){$0$}

\put(6,2){\vector(1,1){1.5}}
\put(18,3.5){\vector(1,-1){1.5}}

\put(8,4.5){$\wedge^{odd}/\wedge_{GL}$}
\put(12,4.5){\vector(1,0){2.5}}
\put(16.5,4.5){$\wedge_{BGL}$}
\put(13,5){\small{$d$}}

\put(6.5,7.5){\vector(1,-1){1.5}}
\put(7.25,7){\small{deR}}
\put(10.5,7){\small{$i$}}
\put(10.5,6){\vector(1,1){1.5}}
\put(14.5,7.5){\vector(1,-1){1.5}}
\put(15.5,7){\small{$ch$}}
\put(17.75,7){\small{deR}}
\put(18.5,6){\vector(1,1){1.5}}

\put(3,8){$H^{odd}(C)$}
\put(12.75,8){$\hat{K}$}
\put(20,8){$H^{even}(C)$}

\put(1.5,10.5){\small{reduction mod $Z$}}
\put(6,9.5){\vector(1,1){1.5}}
\put(10.5,11){\vector(1,-1){1.5}}
\put(11.5,10.5){\small{$j$}}
\put(14,10.5){\small{$\delta$}}
\put(14.0,9.5){\vector(1,1){1.5}}
\put(18,11){\vector(1,-1){1.5}}
\put(19,10.5){\small{$c$}}

\put(8,12){$K(C/Z)$}
\put(12,12){\vector(1,0){2.5}}
\put(16.5,12){$K$}
\put(11.75,12.5){\footnotesize{Bockstein}}

\put(5.5,14.5){\vector(1,-1){1.5}}
\put(18,13){\vector(1,1){1.5}}

\put(4.5,15){$0$}
\put(20,15){$0$}
\end{picture}
\end{center}
where the sequence along the upper boundary may be identified (via $ch
\otimes C$) with the
Bockstein sequence for complex $K$-theory (the long exact sequence
associated to the short exact sequence of coefficients $0 \rightarrow
Z \rightarrow C \rightarrow C/Z \rightarrow 0$), and that along the lower
boundary comes from de Rham theory.  Here, $\wedge_{BGL}$ means all
closed forms cohomologous to Chern characters of complex vector
bundles, and $\wedge_{GL}$ means all closed forms cohomologous to
pull-backs by maps into \mbox{$GL = \textrm{union of the } GL(n,C)$} of the transgression of the Chern
character form.  $\delta$ is the map which simply drops the
connection, and $ch$ is the Chern-Weil map applied to the Chern character polynomial.
The fibre product statement above is related to the commutative
square on the right half of the diagram.

The work's main technical innovation is embodied in Proposition 2.6,
where it is shown that all odd forms modulo $\wedge_{GL}$ arise as the
Chern-Simons difference forms between the trivial connection and
arbitrary connections
on trivial bundles.  A corollary, as implied by the diagram above, is
that every element of $\wedge_{BGL}$ arises after stabilizing as the Chern character form of
some connection in any bundle whose Chern character is the given
cohomology class.  In particular, if a bundle has zero characteristic
classes over $C$, then there is a connection on that bundle,
stabilized by adding in a trivial bundle, with
vanishing Chern-Weil forms.

By considering the simultaneous kernel of $ch$ and $\delta$, the
diagram also shows that the ambiguity in determining a structured
bundle up to stabilizing solely by its characteristic forms and underlying element of
$K$ is measured by a complex torus, the dimension of which is the sum
of the odd Betti numbers of the base manifold.

In showing that the kernel of $ch$ is $K(C/Z)$ we were influenced by
the work of Karoubi [2] and Lott [1], which gave a related description
of $K(C/Z)$ involving a bundle with connection and an extra total odd
form whose $d$ is the Chern character form.  Our proof is based on the
characterization given in Appendix A of the homotopy fibre of a map in
the homotopy category.

We also point out that the
existence of a differential $K$-theory associated to $K$-theory, and
indeed a differential theory associated to any exotic cohomology theory, was
constructed in the paper of Hopkins and Singer [5].  Following their
approach, Freed, as well as Hopkins and Singer and perhaps others like
ourselves, were aware that a model for differential $K$-theory could
be constructed based on pairs $(E,O)$, where $E$ is a bundle with
connection and $O$ is a total odd form with an equivalence relation
generalizing that in [2].  One point of the present work is that this
total odd form may be taken to be zero in the equivalent description of
differential $K$-theory presented here.

There is a word for word variant of the above concerning complex
vector bundles with Hermitian connection.  Now there is a functor
$\hat{K}_{R}$, four natural transformations and the diagram

\vspace{.5cm}
\begin{center}
\setlength{\unitlength}{0.5cm}
\begin{picture}(24,16)\thicklines
\put(5,1){$0$}
\put(20.5,1){$0$}

\put(6,2){\vector(1,1){1.5}}
\put(18,3.5){\vector(1,-1){1.5}}

\put(8,4.5){$\wedge^{odd}_{\wedge_{U}}$}
\put(12,4.5){\vector(1,0){2.5}}
\put(16.5,4.5){$\wedge_{BU}$}
\put(13,5){\small{$d$}}

\put(6.5,7.5){\vector(1,-1){1.5}}
\put(7.25,7){\small{deR}}
\put(10.5,7){\small{$i$}}
\put(10.5,6){\vector(1,1){1.5}}
\put(14.5,7.5){\vector(1,-1){1.5}}
\put(15.5,7){\small{$ch$}}
\put(17.75,7){\small{deR}}
\put(18.5,6){\vector(1,1){1.5}}

\put(3,8){$H^{odd}(R)$}
\put(12.75,8){$\hat{K}_{R}$}
\put(20,8){$H^{even}(R)$}

\put(1.5,10.5){\small{reduction mod $Z$}}
\put(6,9.5){\vector(1,1){1.5}}
\put(10.5,11){\vector(1,-1){1.5}}
\put(11.5,10.5){\small{$j$}}
\put(14,10.5){\small{$\delta$}}
\put(14.0,9.5){\vector(1,1){1.5}}
\put(18,11){\vector(1,-1){1.5}}
\put(19,10.5){\small{$c$}}

\put(8,12){$K(R/Z)$}
\put(12,12){\vector(1,0){2.5}}
\put(16.5,12){$K(Z)$}
\put(11.75,12.5){\footnotesize{Bockstein}}

\put(5.5,14.5){\vector(1,-1){1.5}}
\put(18,13){\vector(1,1){1.5}}

\put(4.5,15){$0$}
\put(20,15){$0$}
\end{picture}
\end{center}

This is discussed briefly in Section 5.  As a corollary, for any
bundle over a closed Riemannian manifold after stabilizing, there is a
unitary connection on the bundle whose Chern-Weil form is the harmonic
representative of the Chern character of the bundle.  Moreover, when
the odd Betti numbers vanish, this structured bundle is unique up to
adding factors with trivial holonomy.

Our model of $\hat{K}$ or $\hat{K}_{R}$  may relate to two questions:

\begin{itemize}
\item[1.]  Up to a natural transformation, are $\hat{K}$ or $\hat{K}_{R}$ uniquely
determined by the diagram, as shown in [7]  in the case of ordinary
differential cohomology?

\item[2.]  Can one enrich the families index theorem by passing from $K$
to $\hat{K}$ or $\hat{K}_{R}$?  c.f. [3], [4], [6].
\end{itemize}

Finally, this model of $\hat{K}$ or $\hat{K}_{R}$ might be helpful for
certain quantum theories and \mbox{$M$-theory}, in which it has already been
observed that actions can be written more appropriately in the
language of differential $K$-theory than in that of differential forms
[6], [8].  In this respect we note Theorem 3.9, showing that $\hat{K}$
and  $\hat{K}_{R}$
satisfy the Mayer-Vietoris property, which relates to locality.

\section*{\S1. Structured Bundles}

Let $[V,\nabla]$ be a complex vector bundle with connection over a
smooth manifold with corners, $X$, and let \mbox{$R \in \wedge^2 (X, \textrm{End}(V))$} denote
its curvature tensor.

Using the Chern-Weil homomorphism, the Chern character of $V$,
$ch(V)$, may be represented by the total complex valued closed form on
$X$, $ch(\nabla)$, defined by 
\[
1.1) \qquad ch(\nabla) = \sum_{j=0} \frac{1}{j!} 
\left( \frac{1}{2 \pi i} \right)^{j} \mathrm{tr}
( \overbrace{R \wedge \cdots \wedge R}^{j} ) \quad \in \wedge^{even}(X,C).
\]
For $t \in [0,1]$ and  $\gamma(t)= \nabla^t$ a smooth curve of
connections, $(\nabla^t)^\prime = A^t \in \wedge^1 (X, \textrm{End}(V))$,
and we set
\[
1.2) \qquad \mathit{cs}(\gamma)= \int_0^1\quad \sum_{j=1} \frac{1}{(j-1)!}
\left( \frac{1}{2 \pi i}\right)^j \mathrm{tr}
( A^{t} \wedge \overbrace{R^{t} \wedge \cdots \wedge R^{t}}^{j-1} ) \quad \in \wedge^{odd}(X,C).
\]
It is a standard fact that
\[
1.3) \qquad \textrm{d} \mathit{cs}(\gamma)= ch(\nabla^1) - ch (\nabla^0).
\]

There is a second formulation of 1.2) which will be useful in what
follows.

Let $\Pi : X \times [0,1] \rightarrow X$ be the standard projection, and
set $W = \Pi^{*}(V)$.  We may construct a connection, $\bar{\nabla}$,
on $W$ by defining ${\bar{\nabla}}_{s} = \nabla^{t}_{\Pi_{*}(s)}$ when
$s$ is tangent to the slice through $t$, and by making
${\bar{\nabla}}_{\partial/\partial t}(\Pi^{*}(f)) = 0$ for $f$ any
cross-section of $V$.

Let $\bar{R}$ be the curvature tensor of $\bar{\nabla}$.  Then, if
$r,s$ are tangent to the slice through $t$,
\[ \begin{array}{llcl}
& \bar{R}_{r,s} & = & R^{t}_{\Pi_{*}(r),\Pi_{*}(s)} \\
1.4) \quad &&& \\
& \bar{R}_{\partial/\partial t, \: s} & = & A^{t}_{\Pi_{*}(s)} . \\
\end{array}  \]

The first is straightforward.  To show the second, let $w \in
W_{(x,t)}$, and extend it to be of the form $\Pi^{*}(f)$, where $f$ is a
cross-section of $V$.  Also extend $s$ to be the lift of a vector
field on $X$.  Clearly $[s,\partial/\partial t] = 0$.  Thus
\[
\bar{R}_{\partial/\partial t, \: s} \, w = \bar{\nabla}_{\partial/\partial t} \bar{\nabla}_{s} w - \bar{\nabla}_{s}
\bar{\nabla}_{\partial/\partial t} w = \bar{\nabla}_{\partial/\partial t}
\bar{\nabla}_{s} w = \frac{d}{dt} \nabla^{t}_{\Pi_{*}(s)} w = A^{t}_{\Pi_{*}(s)} w.
\]
Now, let $\psi_{t} : X \rightarrow X \times [0,1]$ be the slice map,
$\psi_{t}(x) = (x,t)$.  Then by 1.4)
\[ \begin{array}{lcl}
\mathrm{tr} (A^{t} \wedge \overbrace{R^{t} \wedge \cdots \wedge
R^{t}}^{j-1}) & = & \psi^{*}_{t}(\, \mathrm{tr} (i_{\partial/\partial t} \bar{R} \wedge
\overbrace{\bar{R} \wedge \cdots \wedge \bar{R}}^{j-1}) \,) \\
& = & \psi^{*}_{t}(\, i_{\partial/\partial t} (\, \frac{1}{j} \, \mathrm{tr}
(\overbrace{\bar{R} \wedge \cdots \wedge \bar{R}}^{j})) \, ).  \\ 
\end{array} \]
From this we conclude
\[
1.5) \qquad \mathit{cs}(\gamma) = \int_0^1 \psi^{*}_{t} \,
(i_{\partial/\partial t} \, ch(\bar{\nabla})).
\]

The following proposition is almost certainly well known, but we
did not find a reference.

\underline{\textbf{Proposition 1.6}}:  If $\alpha$ and $\gamma$ are
two paths connecting $\nabla^0$ and $\nabla^1$, then
\[ 
\qquad \mathit{cs}(\alpha)= \mathit{cs}(\gamma) + \textrm{exact}. 
\]

\underline{\textbf{Proof}}:  \quad It is sufficient to prove that if
$\gamma$ is a closed path of connections, then $cs(\gamma)$ is exact. \\
By 1.3) $cs(\gamma)$ is obviously closed.  To show it exact we show
that $cs(\gamma)$ integrates to $0$ on every cycle of $X$.

Let $Z$ be such a cycle.  Then by 1.5)
\[
\int_{Z} cs(\gamma) = \int_{Z \times S^{1}} ch(\bar{\nabla}) =
ch(W)(Z \times S^{1}) = \Pi^{*}(ch(V))(Z \times S^{1}) =
ch(V)(\Pi_{*}(Z \times S^{1})) = 0
\]
Thus $cs(\gamma)$ is exact. \hspace{1 cm}  $\blacksquare$

\vspace{0.3cm}

Since $\nabla^0$ and $\nabla^1$ may always be joined by a smooth path,
using Proposition 1.6, we may set
\[
1.7) \qquad \mathit{CS}(\nabla^0, \nabla^1) = cs(\gamma) \, \bmod \mathrm{exact}.
\]
From Proposition 1.6) we also see
\[
1.8) \qquad \mathit{CS}(\nabla^0, \nabla^1) + \mathit{CS}(\nabla^1, \nabla^2) = \mathit{CS}(\nabla^0, \nabla^2).
\]

\underline{\textbf{Definition}}:   $\nabla^0$ and $\nabla^1$ will be called \textbf{equivalent}, and written
$\nabla^0 \sim \nabla^1$, if $\mathit{CS}(\nabla^0, \nabla^1) = 0$. Equation 1.8) shows $\sim$ is an equivalence relation.

\underline{\textbf{Definition}}:  A pair $\mathcal{V} = [V, \{\nabla\}]$, where
$\{\nabla\}$ is an equivalence class of connections on $V$ will be
called a \textbf{structured bundle}.

If $\nabla^W$ is a connection on $W$ and $\sigma : V \rightarrow W$ is
a bundle isomorphism covering the identity map of $X$, $\sigma$
induces $\sigma^{*}(\nabla^W)$, a connection on $V$, and it is easily
seen that $\{\sigma^{*}(\nabla^W)\} = \sigma^{*}(\{\nabla^W\}).$
$\mathcal{V} = [V, \{\nabla^{V}\}]$ and $\mathcal{W} = [W,
\{\nabla^{W}\}]$ are called \textbf{isomorphic} if
$\sigma^{*}(\{\nabla^W\}) = \{\nabla^{V}\}$.

If $\psi:X \rightarrow Y$ is $C^\infty$, and $V$ is a bundle over $Y$
with connections $\nabla^0$ and $\nabla^1$, then, in the usual
manner, $\psi^{*}(\nabla^0)$ and $\psi^{*}(\nabla^1)$ are
connections on $\psi^{*}(V)$.   Clearly
\[
\mathit{CS}(\psi^{*}(\nabla^0), \psi^{*}(\nabla^1)) = \psi^{*}(\mathit{CS}(\nabla^0, \nabla^1)).
\]
Thus, if $\mathcal{V} = [V, \{\nabla\}]$ is a structured bundle over $Y$ then
$\psi^{*}(\mathcal{V}) = [\psi^{*}(V), \{\psi^{*}(\nabla)\}]$ is well
defined as a structured bundle over $X$.

Suppose $\psi_{t} : X \rightarrow Y$ is a smooth 1-parameter family of
maps.  If $\mathcal{V} = [V, \{\nabla\}]$ is a structured bundle over
$Y$, then $\mathcal{V}^{t} = [\psi^{*}_{t}(V),
\psi^{*}_{t}(\{\nabla\})]$ is a 1-parameter family of structured
bundles over $X$.  Assume $t \in [0,1]$ and let $\gamma_{x} : [0,1]
\rightarrow Y$ be the curve $\gamma_{x}(t) = \psi_{t}(x)$.  Let
$\sigma_{t} : \psi^{*}_{0}(V) \rightarrow \psi^{*}_{t}(V)$ be parallel
transport along the curves $\gamma_{t}$.  Then, letting $W =
\psi^{*}_{0}(V)$ and $\nabla^{t} =
\sigma^{*}_{t}(\psi^{*}_{t}(\nabla))$, $\mathcal{W}^{t} = [W, \{\nabla^{t}\}]$
is a 1-parameter family of structured bundles over $X$, isomorphic to
the family $\mathcal{V}^{t}$, having the same underlying vector bundle.

Letting $\gamma^{\prime}_{x}(t)$ denote the tangent vector to
$\gamma_{x}$ at $t$, and using 1.5), we conclude
\[
1.9) \qquad \mathit{CS}(\nabla^{0},\nabla^{1}) = \int_{0}^{1} \psi^{*}_{t}
(i_{\gamma^{\prime}_{x}(t)} \, ch(\nabla)) \, dt.
\]

\vspace{0.3cm}

If $\nabla^V$ and $\nabla^W$ are connections on $V$ and $W$
they determine connections on $V \oplus W$ and $V \otimes
W$, denoted by $\nabla^V \oplus \nabla^W$ and $\nabla^W \otimes
\nabla^W$. For $f,g$ cross-sections in $V$ and $W$, and $r$,
a tangent vector to $X$,
\[ \begin{array}{lcl}
(\nabla^{V} \oplus \nabla^{W})_{r}(f,g) & = & (\nabla^{V}_{r} f,
 \nabla^{W}_{r} g) \\
(\nabla^{V} \otimes \nabla^{W})_{r}(f \otimes g) & = & \nabla^{V}_{r}(f)
 \otimes g + f \otimes \nabla^{W}_{r}(g). 
\end{array} \]
It is well known that
\[  \begin{array}{llcl}
1.10) & ch(\nabla^V \oplus \nabla^W) & = & ch(\nabla^V) + ch(\nabla^W) \\
1.11) & ch(\nabla^V \otimes \nabla^W) & = & ch(\nabla^V) \wedge ch(\nabla^W). 
\end{array} \]

\underline{\textbf{Lemma 1.12}}:  Let $\nabla^{V}, \bar{\nabla}^{V}, \nabla^W,
\bar{\nabla}^{W}$ be connections on the indicated bundles.  Then
\[ \begin{array}{llcl}
a) & \mathit{CS}(\nabla^V \oplus \nabla^W, \bar{\nabla}^V \oplus
\bar{\nabla}^W) & = & \mathit{CS}(\nabla^V,\bar{\nabla}^V) + \mathit{CS}(\nabla^W,\bar{\nabla}^W) \\
b) & \mathit{CS}(\nabla^V \otimes \nabla^W, \bar{\nabla}^V \otimes
\bar{\nabla}^W) & = & ch(\nabla^V) \wedge \mathit{CS}(\nabla^W,\bar{\nabla}^W)
+ ch(\bar{\nabla}^W) \wedge \mathit{CS}(\nabla^V,\bar{\nabla}^V) \\
\end{array} \]

\underline{\textbf{Proof}}:  Using 1.8)
\[
\mathit{CS}(\nabla^{V} \oplus \nabla^W, \bar{\nabla}^V \oplus \bar{\nabla}^W) =
\mathit{CS}(\nabla^{V} \oplus \nabla^W, \nabla^{V} \oplus \bar{\nabla}^W) + \mathit{CS}(\nabla^V \oplus \bar{\nabla}^W, \bar{\nabla}^V \oplus \bar{\nabla}^{W}).
\]
Direct calculation of each term using 1.2) shows a).

Again using 1.8)
\[
\mathit{CS}(\nabla^{V} \otimes \nabla^W, \bar{\nabla}^V \otimes \bar{\nabla}^W)
= \mathit{CS}(\nabla^{V} \otimes \nabla^W, \nabla^{V} \otimes \bar{\nabla}^W) + \mathit{CS}(\nabla^V \otimes \bar{\nabla}^W, \bar{\nabla}^V \otimes \bar{\nabla}^W)
\]
and again from 1.2), direct calculation shows b). \hspace{1 cm} $\blacksquare$

\vspace{0.3cm}

From Lemma 1.12 one immediately sees

\underline{\textbf{Proposition 1.13}}:  If $\mathcal{V} = [V, \{\nabla^{V}\}]$
and $\mathcal{W} = [W, \{\nabla^{W}\}]$ are structured bundles, then the equivalence classes $\{\nabla^V \oplus \nabla^W\}$ and $\{\nabla^V \otimes \nabla^W\}$ are independent of the choices of $\nabla^V \in \{\nabla^{V}\}$ and $\nabla^W \in \{\nabla^{W}\}$, and so
\begin{eqnarray*}
\mathcal{V} \oplus \mathcal{W} & = & [V \oplus W, \{\nabla^V \oplus \nabla^W\}] \qquad \textrm{and} \\
\mathcal{V} \otimes \mathcal{W} & = & [V \otimes W, \{\nabla^V \otimes \nabla^W\}] 
\end{eqnarray*}
are well defined structured bundles.

\underline{\textbf{Definition}}:  We define \textbf{Struct($X$)} to
be the set of isomorphism classes of structured bundles over $X$.  By
Proposition 1.13, the operations $\oplus$ and $\otimes$ make
Struct($X$) an abelian semi-group with commutative, distributive
multiplication.  A smooth map $\psi$ from $X$ to $Y$ induces $\psi^{*}
: \textrm{Struct}(Y) \rightarrow \textrm{Struct}(X)$ preserving these operations.  Thus,
Struct is a \textbf{functor} on the category of smooth compact
manifolds with corners into that of commutative semi-rings.

We conclude from 1.3) that $ch : \textrm{Struct}(X) \rightarrow \wedge^{even}(X)$ is a well defined natural transformation, and from 1.10) and 1.11)
\[
\begin{array}{lcll}
& ch(\mathcal{V} \oplus \mathcal{W}) & = & ch(\mathcal{V}) +
ch(\mathcal{W}) \\
1.14) &&& \\
& ch(\mathcal{V} \otimes \mathcal{W}) & = & ch(\mathcal{V}) \wedge
ch(\mathcal{W}). \\
\end{array}
\]

\underline{\textbf{Definition}}:  A connection $\nabla$ on $V$ will
be called \textbf{flat} if its holonomy around every closed path is
the identity.  This implies the curvature $R \equiv 0$ and that $V$ is
 naturally isomorphic to the product bundle with the trivial product
connection.  $\mathcal{V} = [V, \{\nabla\}]$ will be called flat if some ${\nabla} \in \{\nabla\}$ is flat. Since any two such of dim $n$ are isomorphic, we shall denote this isomorphism class by $[n] \in \textrm{Struct}(X)$.

\vspace{0.3cm}

The following theorem is based on a related result in [11], stated
without giving the proof.  We employ that proof here in Lemma 1.16 below.

\underline{\textbf{Theorem 1.15}}:  Given any $\mathcal{V} \in
\textrm{Struct}(X)$ there is a $\mathcal{W} \in \textrm{Struct}(X)$
such that $\mathcal{V} \oplus \mathcal{W} = [n]$ for some $n$.  Any such $\mathcal{W}$ will be called an \textbf{inverse} of $\mathcal{V}$.

To prove the Theorem we need

\underline{\textbf{Lemma 1.16}}:  Let $\nabla$ be a connection on $V
\oplus W$ with curvature $R$.  Let $\nabla^{V}$ and $\nabla^W$ be the
connections on $V$ and $W$ induced by $\nabla$. E.g. if $\Pi^{V} :
V \oplus W \rightarrow V$ is the projection, and $f$ is a cross-section
in $V$ then $\nabla^{V}_{r} f = \Pi^{V}(\nabla_{r} f)$.  Suppose
$R_{r,s}(V) \subseteq V$ and $R_{r,s}(W) \subseteq W$ for all
tangent vectors $r,s$ at any point of $X$.  Then,
\[
\nabla^{V} \oplus \nabla^W \sim \nabla.
\]

\underline{\textbf{Proof}}:  We may write
\[
\nabla= \nabla^{V} \oplus \nabla^{W} + A
\]
where $A \in \wedge^{1}(X, \textrm{End}(V \oplus W))$.  For $f$ a cross-section in $V$ we see
\[
A_{r} f = \nabla_{r} f - \Pi^{V}(\nabla_{r} f) = \Pi^{W}(\nabla_{r} f) \in W.
\]
As the same holds for $W$, we see
\[
1.17) \qquad A_{r}(V) \subseteq W \qquad \textrm{and} \qquad A_{r}(W) \subseteq V.
\]
Setting $\bar{\nabla} = \nabla^{V} \oplus \nabla^W$, let $\bar{R}$
denote its curvature and $\bar{d}$ denote its exterior differentiation
operator.  Since $\bar{\nabla}$ preserves $V$ and $W$, 1.17) implies
\[
1.18) \qquad \bar{d}A_{r,s}(V) \subseteq W \qquad and \qquad
\bar{d}A_{r,s}(W) \subseteq V.
\]
The usual formula computing the curvature of one connection from that of another shows
\[
R= \bar{R} + A \wedge A + \bar{d} A.
\]
By hypothesis, $R$ preserves $V$ and $W$. So does $\bar{R}$, being the
curvature of a direct sum connection, and so does $A \wedge A$ by
1.17).  This implies that $\bar{d} A$ preserves them as well, but 1.18)
shows the opposite. Thus $\bar{d}A = 0$ and
\[
1.19) \qquad R = \bar{R} + A \wedge A.
\]

Let $\nabla^{t} = \bar{\nabla} + t A$, a curve of connections joining
$\nabla^V \oplus \nabla^W$ to $\nabla$.  Letting $R^{t}$ denote the associated curvature, we see from 1.19)
\[
1.20) \qquad R^{t}= \bar{R} + t^{2} A \wedge A.
\]
In the notation of 1.2), $A^{t} = (\nabla^{t})^\prime = A$, and so the $\mathit{CS}$ integrand consists of terms of the form
\[
\mathrm{tr} (A \wedge \overbrace{R^{t} \wedge \cdots \wedge R^{t}}^{j-1}).
\]
But, by 1.20) $R^{t}$ preserves both $V$ and $W$, and, since $A$
reverses them, all such trace terms must vanish. Thus $\mathit{CS}(\nabla,
\bar{\nabla}) = 0$.  \hspace{1 cm} $\blacksquare$

\vspace{0.3cm}

\underline{\textbf{Proof of Theorem 1.15}}:  

The classifying spaces $B_k GL (n,C) = GL (n+k, C) / GL (n, C) \times GL (k, C)$ come with
natural bundles, $V^n$ and $V^k$, of dimension $n$ and $k$, and
connections $\nabla^n$ and $\nabla^k$ induced by the standard flat
connection on $V^n \oplus V^k$.  Lemma 1.16 shows that $\mathcal{V}^n= [V^n, \{\nabla^n\}]$ and $\mathcal{V}^k=[V^k,\{\nabla^k\}]$ are inverses of each other.

The theorem of Narasimhan-Ramanan [9] shows that for sufficiently
large $k$, an $n$-dim $\mathcal{V} \in \textrm{Struct}(X)$ may be
obtained as the pull-back of $\mathcal{V}^n$ via a $C^\infty$ map of
$X \rightarrow B_k GL (n,C)$. The pull-back of $\mathcal{W}^k$
will then be an inverse of $\mathcal{V}$ in the sense of Theorem
1.15.  \hspace{1 cm} $\blacksquare$

\section*{\S2. The Stably Trivial Case}

Let $GL = \underset{n}{\lim} \, \, GL (n,C)$, the stabilized complex general linear group and $\mathcal{G}$ its Lie algebra. $\mathcal{G}$ consists of complex valued matrices, all but a finite number of whose entries are $0$. Let $\theta \in \wedge^{1} (GL,\mathcal{G})$ denote the canonical left invariant $\mathcal{G}$-valued form on $GL$. Set
\[
2.1) \qquad \Theta = \sum_{j=1} b_{j} \, \mathrm{tr} (\overbrace{\theta
\wedge \theta \wedge \cdots \wedge \theta}^{2j-1}) \, \in \wedge^{odd}(GL)
\]
where
\[ 
b_{j} = \frac{1}{(j-1)!} \left(\frac{1}{2 \pi i}\right)^{j} \, \int_0^1 (t^{2} - t)^{j-1} \, dt.
\]

It is well known that $\Theta$ is a bi-invariant closed odd form, and
the free abelian group generated by all distinct products of its
components represent the entire complex cohomology ring of $GL$. We
define $\wedge_{GL} \subseteq \wedge^{odd}$ by
\[
2.2) \quad \wedge_{GL}(X) = \{g^{*} (\Theta)\} + \wedge_{exact}^{odd}
\]
where $g : X \rightarrow GL$ runs through all $C^\infty$ maps.

Note that if $g,h$ map $X$ into $GL$, then $g \, \oplus \, h : X \rightarrow
GL$ may be defined, and $(g \, \oplus \, h)^{*}(\Theta) = g^{*}(\Theta) +
h^{*}(\Theta)$.  Moreover, $(g^{-1})^{*}(\Theta) = -g^{*}(\Theta)$.
Thus $\wedge_{GL}(X)$ is an abelian group.

\vspace{0.3cm}

\underline{\textbf{Lemma 2.3}}:  Let $V$ be a trivial bundle with the
two flat connections $\nabla$ and $\bar{\nabla}$.  Then
\[
\mathit{CS}(\nabla, \bar{\nabla}) \in \wedge_{GL} / \wedge_{exact}^{odd}.
\]

\underline{\textbf{Proof}}:  Since $\nabla$ and $\bar{\nabla}$ each
have trivial holonomy, we can find a cross-section $g \in \mathrm{Aut}(V)$ such that
\[
\bar{\nabla}_{t}(f) = g^{-1} (\nabla_{t} (g(f))).
\]
Expressing $g$ as a matrix with respect to a $\nabla$-parallel framing of $V$, we see
\[
\bar{\nabla} = \nabla + g^{-1} \, dg.
\]
Now, regarding $g : X \rightarrow GL$, one easily sees that $g^{-1} d g
= g^{*}(\theta)$.  Thus
\[
\bar{\nabla} = \nabla + g^{*}(\theta).
\]
Setting  $\bar{\nabla}^{t} = \nabla + t g^{*} (\theta)$, we see
\[
\bar{R}^{t} = R + t \, d g^{*}(\theta) + t^{2} \, g^{*}(\theta) \wedge g^{*}(\theta).
\]
But, either calculating on $GL$, or directly with $g^{-1} dg$, we see
that $ dg^{*}(\theta) = -g^{*}(\theta) \wedge g^{*}(\theta)$.
Moreover, since $\nabla$ has trivial holonomy, $R \equiv 0$. Thus
\[
\bar{R}^{t} = (t^{2} - t) g^{*}(\theta) \wedge g^{*}(\theta).
\]
It then follows from 1.2) that $\mathit{CS}(\nabla,\bar{\nabla}) =
g^{*}(\Theta)$.  \hspace{1 cm} $\blacksquare$

\vspace{0.3cm}

\underline{\textbf{Definition}}:
\[
\mathrm{Struct}_{\mathrm{ST}}(X) = \{ \, [V, \{\nabla\}] \in
\mathrm{Struct}(X) \mid \textrm{V is stably trivial} \, \}.
\]
For $\mathcal{V} \in \mathrm{Struct}_{\mathrm{ST}}(X)$, let $F$ and
$H$ be trivial bundles such that $V \oplus F = H$ and let $\nabla^{F},
\nabla^{H}$ be flat connections on $F$ and $H$.  We define
\[
\widehat{\mathit{CS}} : \textrm{Struct}_{\textrm{ST}}(X) \rightarrow
\wedge^{odd} / \wedge_{GL}
\]
by
\[
\widehat{\mathit{CS}}(\mathcal{V}) = \mathit{CS}(\nabla^{H}, \nabla \oplus
\nabla^{F}) \, \bmod \wedge_{GL} / \wedge_{exact}^{odd}.
\]

\vspace{0.3cm}

\underline{\textbf{Proposition 2.4}}: $\widehat{\mathit{CS}}$ is a well defined homomorphism.

\vspace{0.3cm}

\underline{\textbf{Proof}}:  Suppose $\bar{F}, \bar{H},
\nabla^{\bar{F}}, \nabla^{\bar{H}}$ are another pair of trivial
bundles with flat connections with \mbox{$V \oplus \bar{F}= \bar{H}$}.
Using 1.7), Lemma 1.12 and Lemma 2.3, and working 
$\bmod \, \wedge_{GL}$, we see
\begin{eqnarray*}
\mathit{CS}(\nabla^{\bar{H}}, \nabla \oplus \nabla^{\bar{F}}) & = &
\mathit{CS}(\nabla^{\bar{H}} \oplus \nabla^{F}, \nabla \oplus
\nabla^{\bar{F}} \oplus \nabla^{F}) \\
& = & \mathit{CS}(\nabla^H \oplus \nabla^{\bar{F}}, \nabla \oplus
\nabla^{\bar{F}} \oplus \nabla^F) \\
& = & \mathit{CS}(\nabla^H, \nabla \oplus \nabla^F). 
\end{eqnarray*}
Thus $\widehat{\mathit{CS}}$ is well defined. That $\widehat{\mathit{CS}}$ is a homomorphism follows
immediately from Lemma 1.10.  $\blacksquare$

\vspace{0.3cm}

\underline{\textbf{Definition}}:   $\mathcal{V} \in
\mathrm{Struct}(X)$ is called \textbf{stably flat} if there exists
flat $\mathcal{F}$ and $\mathcal{H}$ such that $\mathcal{V} \oplus
\mathcal{F}= \mathcal{H}$.  The set of these objects will be denoted by
$\textrm{Struct}_{\textrm{SF}}(X)$.  Clearly
$\textrm{Struct}_{\textrm{SF}}(X) \subseteq \textrm{Struct}_{\textrm{ST}}(X)$ and is a sub semi-group.

\vspace{0.3cm}

\underline{\textbf{Proposition 2.5}}:  $\textrm{ker}(\widehat{\mathit{CS}}) = \textrm{Struct}_{\textrm{SF}}(X)$.

\vspace{0.3cm}

\underline{\textbf{Proof}}:  Obviously $\textrm{Struct}_{\textrm{SF}}
\subseteq \textrm{ker}(\widehat{\mathit{CS}})$.  Now suppose $\widehat{\mathit{CS}}(\mathcal{V})=
0$. Let $F, H, \nabla^{F}$ and $\nabla^{H}$ be as in the definition of
$\widehat{\mathit{CS}}$.  Now, $\widehat{\mathit{CS}}(\mathcal{V}) = 0$ implies
\[
\mathit{CS}(\nabla^{H}, \nabla \oplus \nabla^{F}) = g^{*}(\Theta) \bmod \wedge_{exact}^{odd}
\]
for some $g : X \rightarrow GL$.  Again as in the proof of Lemma 2.3,
choosing a $\nabla^{H}$-parallel framing of $H$, we may regard $g \in
\mathrm{Aut} (H)$ and set 
\[
{\bar{\nabla}}^{H} = g^{-1} (\nabla^{H} \circ g).
\]
As in the Lemma we see $\mathit{CS}(\nabla^{H}, {\bar{\nabla}}^{H}) =
g^{*}(\Theta)$ and thus $\mathit{CS}({\bar{\nabla}}^{H}, \nabla^{H})
= -g^{*}(\Theta)$.  Therefore
\begin{eqnarray*}
\mathit{CS}({\bar{\nabla}}^{H}, \nabla \oplus \nabla^{F}) & = &
\mathit{CS}({\bar{\nabla}}^{H}, \nabla^H) +
\mathit{CS}(\nabla^{H}, \nabla \oplus \nabla^{F}) \\
&= & -g^{*}(\Theta)+ g^{*}(\Theta) = 0 \, \bmod \mathrm{exact}. 
\end{eqnarray*}
Setting $\bar{\mathcal{H}} = [H, \{ {\bar{\nabla}}^{H} \}]$ and
$\mathcal{F} = [F, \{ \nabla^{F} \}]$ we see $\mathcal{V} \oplus
\mathcal{F} = \bar{\mathcal{H}}$ and thus $\mathcal{V} \in
\textrm{Struct}_{\textrm{SF}}(X)$. \hspace{1 cm} $\blacksquare$

\vspace{0.3cm}

\underline{\textbf{Proposition 2.6}}:  $\textrm{Im}(\widehat{\mathit{CS}}) = \wedge^{odd}(X)/\wedge_{GL}(X)$.

\vspace{0.3cm}

\underline{\textbf{Proof}}:  If $L$ is a trivialized line bundle over
$X$ then any connection on $L$ is simply a complex valued 1-form, $w$.
Since $w \wedge w = 0$, the associated curvature, $R^{w}$, is $dw$, and
$\{w\} = \{w + df \, | \, f \in C^{\infty}(X,C)\}$.

Let ${\mathcal{L}}_{w} = [L, \{w\}]$.  Using $tw$ as a curve of
connections joining $w$ to the trivial connection, and noting that
$R^{tw} = t R^{w}$,  1.2) shows
\[
2.7) \qquad  \widehat{\mathit{CS}}({\mathcal{L}}_{w}) =  \sum_{j=1} \frac{1}{j!} 
\left( \frac{1}{2 \pi i} \right)^{j} w \wedge (dw)^{j-1}.
\]
We first suppose $X = R^n$.  If $w = f \, dx$ then $w \wedge dw = 0$
and thus $\widehat{\mathit{CS}}({\mathcal{L}}_{fdx}) = f dx$.  Moreover, since $\widehat{\mathit{CS}}$ is a homomorphism
\[
\widehat{\mathit{CS}} \left(\sum_{i} \oplus {\mathcal{L}}_{f_{i} dx_{i}}\right) = \sum f_{i} d x_{i}.
\]
Thus $\wedge^{1}(R^n)/ \wedge_{G}(R^n) \subseteq \textrm{Im}(\widehat{\mathit{CS}})$.

Proceeding by induction on $k$, suppose
\[
2.8) \qquad \left( \sum_{j=1}^{k} \wedge^{2j-1}(R^{n}) \right) /
\wedge_{GL}(R^n) \, \subseteq \textrm{Im}(\widehat{\mathit{CS}}).
\]
Let $w = x_1 dx_2 + x_3 dx_4 + \cdots + x_{2k-1} dx_{2k} + f dx_{2k+1}$.

\underline{\textbf{Claim}}:  $w \wedge (dw)^k = (k+1)! f  dx_1 \wedge
\cdots \wedge dx_{2k+1} + \textrm{exact}$.

To show this, let $\gamma = dx_1 \wedge dx_2 + \cdots + dx_{2k-1}
\wedge dx_{2k}$, \quad and note
\[
dw = \gamma + df \wedge dx_{2k+1} \qquad \Rightarrow \qquad (dw)^k = (\gamma + df \wedge dx_{2k+1})^k. 
\]
Since all powers of $df \wedge dx_{2k}$ vanish,
\begin{eqnarray*}
(dw)^k & = & \gamma^k + k \gamma^{k-1} \wedge df \wedge dx_{2k+1} = k!
 dx_1 \wedge \cdots \wedge dx_{2k} + \\
& & k! \, \left[ \sum_{j=1}^{k} dx_1 \wedge dx_2 \wedge
\cdots \wedge \widehat{dx_{2j-1} \wedge dx_{2j}} \wedge \cdots \wedge
 dx_{2k-1} \wedge dx_{2k} \right]
\wedge df \wedge dx_{2k+1} 
\end{eqnarray*}
Thus,
\begin{eqnarray*}
w \wedge (dw)^{k} & = & k! f dx_1 \wedge \cdots \wedge dx_{2k+1} + \\
& & k! \, \left[ \sum_{j=1}^{k} dx_{1} \wedge \cdots \wedge dx_{2j-2}
\wedge x_{2j-1} \wedge dx_{2j} \wedge \cdots \wedge dx_{2k} \right]
\wedge df \wedge dx_{2k+1} \\
& = & k! \, f dx_1 \wedge \cdots \wedge dx_{2k+1} -
 k! \, \sum_{j=1}^{k} dx_1 \wedge \cdots \wedge dx_{2j-1} \wedge
x_{2j-1} df \wedge dx_{2j} \wedge \cdots \wedge dx_{2k+1} \\
& = & (k+1)! \, f dx_1 \wedge \cdots \wedge dx_{2k+1} + \textrm{exact}. 
\end{eqnarray*}
Thus, working $\bmod \wedge_{GL}(R_n)$,
\vspace{.25cm}
\[
\widehat{\mathit{CS}}(\mathcal{L}_{(2 \pi i)^{k+1}w}) = f dx_1 \wedge \cdots \wedge dx_{2k+1} + \theta,
\]
where 
\[
\theta \in \sum_{j=1}^{k} \wedge^{2j-1}(R^{n}).
\]
By induction, $\theta=\widehat{\mathit{CS}} (\mathcal{V})$ for some $\mathcal{V} \in
\textrm{Struct}_{\textrm{ST}}(R^{n})$.  Theorem 1.15 shows
$\mathcal{V}$ has an inverse ${\mathcal{V}}^{-1}$.  Clearly
${\mathcal{V}}^{-1} \in \textrm{Struct}_{\textrm{ST}}(R^n)$ and by
Proposition 2.4, $\widehat{\mathit{CS}}({\mathcal{V}}^{-1}) = -\theta$. Thus
\[
\widehat{\mathit{CS}}(\mathcal{L}_{(2\pi i)^{k+1}w} \oplus \mathcal{V}^{-1}) = f dx_1
\wedge \cdots \wedge dx_{2k+1}.
\]
The general element of $\wedge^{2k+1}(R^n)$ is the sum of such terms, and thus is the image under $\widehat{\mathit{CS}}$ of the direct sum of the inverse images of each of these terms.

For the general case let $ \psi : X \rightarrow R^n$ be an imbedding.
Since $\psi^{*} : \wedge^{odd}(R^n) \rightarrow \wedge^{odd}(X)$ is
onto, and $\psi^{*}(\wedge_{GL}(R^n)) \subseteq
\wedge_{GL}(X)$, $\psi^{*} : \wedge^{odd}(R^n) /
\wedge_{GL}(R_n) \rightarrow \wedge^{odd}(X) /
\wedge_{GL}(X)$ is onto.  Moreover,
$\psi^{*}(\textrm{Struct}_{\textrm{ST}}(R^n)) \subseteq
\textrm{Struct}_{\textrm{ST}}(X)$, and finally $\widehat{\mathit{CS}} \circ \psi^{*} =
\psi^{*} \circ \widehat{\mathit{CS}}$.  Thus if $\rho \in \wedge^{odd}(X) /
\wedge_{GL}(X)$, we can find $\bar{\rho} \in \wedge^{odd}(R_n) /
\wedge_{GL}(R^n)$ with $\psi^{*}(\bar{\rho}) = \rho$.  By
the special case, $\bar{\rho} = \widehat{\mathit{CS}} (\mathcal{V})$ for some $\mathcal{V} \in \textrm{Struct}_{\textrm{ST}}(R^n)$. Then
\[
\rho = \psi^{*} (\widehat{\mathit{CS}} (\mathcal{V})) =
\widehat{\mathit{CS}}(\psi^{*}(\mathcal{V})). \qquad \blacksquare
\]

From Propositions 2.4, 2.5, 2.6 we see

\underline{\textbf{Theorem 2.7}}:
\[
\widehat{\mathit{CS}} : \textrm{Struct}_{\textrm{ST}}(X) /
\textrm{Struct}_{\textrm{SF}}(X) \quad \underrightarrow{\cong} \quad \wedge^{odd}(X) / \wedge_{GL}(X).
\]

\section*{\S3. \, $\hat{K}(X)$}

Using the standard construction of $K$, which transforms an abelian
semi-group into a group, we define
\[
\hat{K} = K(\textrm{Struct}(X)).
\]
$\hat{K}(X)$ is the free abelian group generated by isomorphism
classes of structured bundles, modulo the relation $\mathcal{V} +
\mathcal{W} - (\mathcal{V} \oplus \mathcal{W})$.  Equivalently
defined, $\hat{K}(X)$ is the quotient of the semi-group under $\oplus$
consisting of all pairs $(\mathcal{V}, \mathcal{W})$ modulo the sub
semi-group consisting of pairs $(\mathcal{V}, \mathcal{V})$.  Since $(0, \mathcal{V})$
is obviously the additive inverse of $(\mathcal{V}, 0)$, we write
$(\mathcal{V}, \mathcal{W})$ as $\mathcal{V} - \mathcal{W}$.

\pagebreak

Using Theorem 1.15 it is straightforward using the pairs definition to show

\begin{tabular}{ll}
3.1) & Every element of $\hat{K}(X)$ is of the form $\mathcal{V} -
[n]$. \\
& \\
3.2) & $\mathcal{V} -[n] = 0 \Leftrightarrow \mathcal{V}$ is stably
flat and $n = \textrm{dim}(\mathcal{V})$. \\
\end{tabular}

Again using the pairs definition, one sees that $\otimes$ is well defined in $\hat{K}(X)$, and thus
$\hat{K}(X)$ becomes a commutative ring.  (Defining $(\mathcal{V},\mathcal{W}) \otimes
(\mathcal{V}^{\prime},\mathcal{W}^{\prime})$ to be $(\mathcal{V} \otimes \mathcal{V}^{\prime} \oplus \mathcal{W} \otimes
\mathcal{W}^{\prime}, \mathcal{W} \otimes \mathcal{V}^{\prime} \oplus \mathcal{V} \otimes \mathcal{W}^{\prime})$ one sees
$\{(\mathcal{W},\mathcal{W})\}$ is an ideal.)

We define $\wedge_{BGL} \subseteq \wedge^{even}$ by
\[
\wedge_{BGL}(X) = \{ ch(\mathcal{V})\} + \wedge^{even}_{exact}
\]
where $\mathcal{V}$ ranges over all elements of $\textrm{Struct}(X)$.
From 1.10) and 1.11) and Theorem 1.15 we see that $\wedge_{BGL}(X)$ is a
commutative ring.

By analogy with the definition of $\wedge_{GL}$, and using the theorem
of Narasimhan-Ramanan [9], we could alternatively have defined
\[
\wedge_{BGL}(X) = \{\phi^{*}(\Omega)\} + \wedge^{even}_{exact}
\]
where $\phi : X \rightarrow BGL$ ranges over all $C^{\infty}$ maps, and
$\Omega$ is the Chern character form of the standard connection on the
classifying bundle over $BGL$.

Clearly, $ch$ extends to $\hat{K}(X)$, and maps it to
$\wedge_{BGL}(C)$.  We also define
\[
\delta : \hat{K}(X) \rightarrow K(X)
\]
by
\[
\delta( \, [V,\{\nabla\}] - [W,\{\bar{\nabla}\}] \, ) = V - W.
\]
Letting $c : K(X) \rightarrow H^{even}(X,C)$ be the natural
transformation defined by the Chern character, and $\textrm{deR} :
\wedge_{BGL}(X) \rightarrow H^{even}(X,C)$ be that defined by the
de Rham Theorem, we see

\vspace{.5cm}
\begin{center}
\setlength{\unitlength}{0.5cm}
\begin{picture}(30,7)\thicklines
\put(8.25,0.5){$\wedge_{BGL}(X)$}
\put(6.5,3){\vector(1,-1){1.5}}
\put(7.5,2.5){\small{$ch$}}
\put(10.5,2.5){\small{deR}}
\put(11,1.5){\vector(1,1){1.5}}
\put(0,3.5){3.3)}
\put(4,3.5){$\hat{K}(X)$}
\put(9,3.5){$\otimes$}
\put(13,3.5){$H^{even}(X,C)$}
\put(6.5,5.5){\small{$\delta$}}
\put(6.5,4.5){\vector(1,1){1.5}}
\put(11,6){\vector(1,-1){1.5}}
\put(12,5.5){\small{c}}
\put(8.5,6.5){$K(X)$}
\end{picture}
\end{center}
is a commutative diagram.

\vspace{0.3cm}

\underline{\textbf{Proposition 3.4}}:  \quad $\textrm{ker}(\delta) \,
\cong \,
\wedge^{odd}(X)/\wedge_{GL}(X)$.

\vspace{0.3cm}

\underline{\textbf{Proof}}:  Define $\Gamma :
\textrm{Struct}_{ST}/\textrm{Struct}_{SF} \rightarrow \hat{K}$ by
\[
\Gamma(\{\mathcal{V}\}) = \mathcal{V} - \textrm{dim}(\mathcal{V}).
\]
By 3.2), $\Gamma$ is well defined and is an injection.  Moreover, it
is clear that $\textrm{Im}(\Gamma) = \textrm{ker}(\delta)$.  Thus from
Theorem 2.7,
\[
\Gamma \circ \widehat{\mathit{CS}}^{-1} : \wedge^{odd}(X)/\wedge_{GL}(X) \,
\stackrel{\cong}{\longrightarrow} \, \textrm{ker}(\delta).  \qquad \blacksquare
\]

\vspace{0.3cm}

Let $i = \Gamma \circ \widehat{\mathit{CS}}^{-1}$.  Since $\delta$ is clearly onto,
\[
3.5) \qquad 0 \longrightarrow \wedge^{odd}(X)/\wedge_{GL}(X)
\stackrel{i}{\longrightarrow} \hat{K}(X)
\stackrel{\delta}{\longrightarrow} K(X) \longrightarrow 0
\]
is an exact sequence.

\vspace{0.3cm}

\underline{\textbf{Proposition 3.6}}:  \quad $ch \circ i = d$, and
$ch$ is onto.

\vspace{0.3cm}

\underline{\textbf{Proof}}: \quad To show the first, note that from
the definition of $\widehat{\mathit{CS}}$,
\[
d\widehat{\mathit{CS}}(\mathcal{V}) = ch(\mathcal{V}) - \textrm{dim}(\mathcal{V}) = ch(\Gamma(\mathcal{V}))
\]
for any $\mathcal{V} \in \textrm{Struct}_{ST}(X)$.  Thus, for $\theta
\in \wedge^{odd}(X)$ and $\{\theta\}$ its equivalence class mod
$\wedge_{GL}(X)$,
\[
ch( i(\{\theta\}) ) = ch( \Gamma(\widehat{\mathit{CS}}^{-1}(\{\theta\})) ) =
d( \{\theta\} ) = d\theta.
\]

To show the second, let $\mu \in \wedge_{BGL}(X)$.  By definition,
$\exists \, \, \mathcal{V} \in \textrm{Struct}(X)$ and $\theta \in
\wedge^{odd}$ so that $\mu = ch(\mathcal{V}) + d\theta$.  By the
above, $\mu = ch(\mathcal{V} + i(\{\theta\})).  \qquad \blacksquare$

\vspace{0.3cm}

Let $\textrm{deR} : H^{odd}(X,C) \rightarrow
\wedge^{odd}(X)/\wedge_{GL}(X)$ be the obvious map induced by the de
Rham Theorem.  Since the image of deR consists of closed forms,
$d \circ \textrm{deR} = 0$, which by Proposition 3.6, implies $ch
\circ i \circ \textrm{deR} = 0$.  Thus, $i \circ
\textrm{deR}(H^{odd}(X,C)) \subseteq \textrm{ker}(ch)$.  We have now
established

\underline{\textbf{Proposition 3.7}}:  The following diagram of
functors and natural transformations is commutative, and its diagonals
are exact.

\vspace{.5cm}
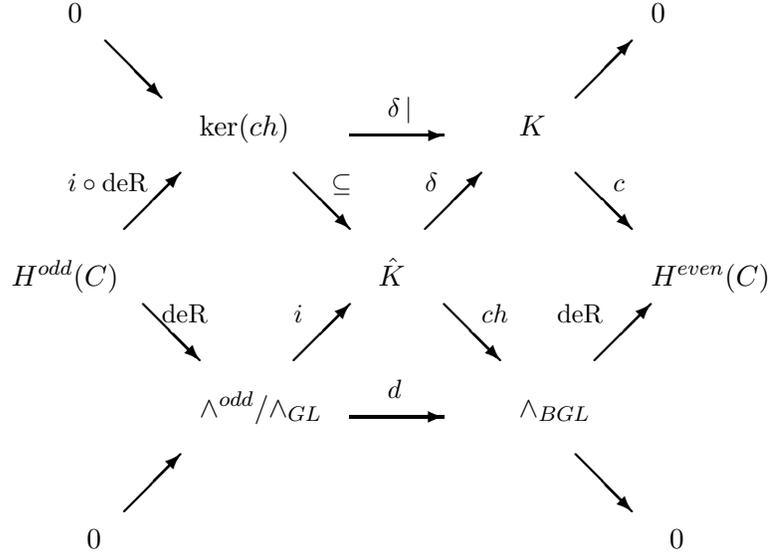
\begin{center}
\setlength{\unitlength}{0.5cm}
\begin{picture}(24,16)\thicklines
\put(5,1){$0$}
\put(20.5,1){$0$}

\put(6,2){\vector(1,1){1.5}}
\put(18,3.5){\vector(1,-1){1.5}}

\put(8,4.5){$\wedge^{odd}/\wedge_{GL}$}
\put(12,4.5){\vector(1,0){2.5}}
\put(16.5,4.5){$\wedge_{BGL}$}
\put(13,5){\small{$d$}}

\put(6.5,7.5){\vector(1,-1){1.5}}
\put(7,7){\small{deR}}
\put(10.5,7){\small{$i$}}
\put(10.5,6){\vector(1,1){1.5}}
\put(14.5,7.5){\vector(1,-1){1.5}}
\put(15.5,7){\small{$ch$}}
\put(17.5,7){\small{deR}}
\put(18.5,6){\vector(1,1){1.5}}

\put(3,8){$H^{odd}(C)$}
\put(12.75,8){$\hat{K}$}
\put(20,8){$H^{even}(C)$}

\put(4.5,10.5){\small{$i \circ \textrm{deR}$}}
\put(6,9.5){\vector(1,1){1.5}}
\put(10.5,11){\vector(1,-1){1.5}}
\put(11.5,10.5){\small{$\subseteq$}}
\put(14,10.5){\small{$\delta$}}
\put(14.0,9.5){\vector(1,1){1.5}}
\put(18,11){\vector(1,-1){1.5}}
\put(19,10.5){\small{$c$}}

\put(8,12){$\textrm{ker}(ch)$}
\put(12,12){\vector(1,0){2.5}}
\put(16.5,12){$K$}
\put(13,12.5){\small{$\delta \, |$}}

\put(5.5,14.5){\vector(1,-1){1.5}}
\put(18,13){\vector(1,1){1.5}}

\put(4.5,15){$0$}
\put(20,15){$0$}
\end{picture}
\end{center}

\underline{\textbf{Corollary 3.8}}:  The outside sequences
\[
\begin{array}{ccccccc}
H^{odd}(C) & \stackrel{i \circ \textrm{deR}}{\longrightarrow} &
\textrm{ker}(ch) & \stackrel{\delta \, |}{\longrightarrow} & K &
\stackrel{c}{\longrightarrow} & H^{even}(C)\\
\\
H^{odd}(C) & \stackrel{\textrm{deR}}{\longrightarrow} &
\wedge^{odd}/\wedge_{GL} & \stackrel{d}{\longrightarrow} & \wedge_{BGL} &
\stackrel{\textrm{deR}}{\longrightarrow} & H^{even}(C)\\
\end{array}
\]
are exact.

\vspace{0.3cm}

\underline{\textbf{Proof}}: \quad Exactness of the first follows from
diagram chasing, and that of the second from the \\ 
de Rham Theorem.
$\qquad \blacksquare$

\vspace{0.3cm}

In the decomposition below we assume that $D$ is a codimension zero or one
submanifold with collar neighborhoods in each of $A$ and $B$.  Thus a smooth form on $D$ can be extended to a
smooth form on either $A$ or $B$.

\underline{\textbf{Theorem 3.9 (Mayer-Vietoris)}}:

Let $A,B \subseteq X$ with $A \cap B = D$ and $A \cup B = X$.  If
$\mu_{A} \in \hat{K}(A)$ and $\mu_{B} \in \hat{K}(B)$ with $\mu_{A} | D
= \mu_{B} | D$, then there exists $\mu \in \hat{K}(X)$ with $\mu | A =
\mu_{A}$ and $\mu | B = \mu_{B}$.

\vspace{0.3cm}

\underline{\textbf{Proof}}: \quad  Following the diagram in
Proposition 3.7, since $ \delta(\mu_{A}) \, | \, D = \delta
(\mu_{B}) \, | \, D$, the Mayer-Vietoris property for $K$ produces $V - [n]
\in K(X)$ with $(V - [n]) \, | \, A = \delta(\mu_{A})$ and $(V - [n]) \,
| \, B = \delta(\mu_{B})$.  Choose $\bar{\mu} \in \hat{K}(X)$ with
$\delta(\bar{\mu}) = V - [n]$.

Now, $\delta(\bar{\mu} \, | \, A) = \delta(\bar{\mu}) \, | \, A =
\delta(\mu_{A})$, and similarly for $B$.  Thus, by the diagram
\[ \begin{array}{llcl}
& \bar{\mu} \, | \, A & = & \mu_{A} + i(\{\alpha_{A}\}) \\
\ast) \quad &&& \\
& \bar{\mu} \, | \, B & = & \mu_{B} + i(\{\alpha_{B}\}) \\
\end{array}  \]

where $\alpha_{A}, \, \alpha_{B} \in \wedge^{odd}(A), \, \wedge^{odd}(B)$
and $\{\alpha_{A}\}$, $\{\alpha_{B}\}$ represent their equivalence
classes mod $\wedge_{GL}(A)$,  $\wedge_{GL}(B)$.

By the above,
\[
i(\{\alpha_{A} \, | \, D\}) - i(\{\alpha_{B} \, | \, D\}) =
i(\{\alpha_{A}\}) \, | \, D - i(\{\alpha_{B}\}) \, | \, D = (\bar{\mu}
\, | \, A) \, | \, D - (\bar{\mu} \, | \, B) \, | \, D - \mu_{A} \, | \,
D + \mu_{B} \, | \, D.
\]
The first pair vanishes since each term is $\bar{\mu} \, | \, D$, and
the second pair vanishes by hypothesis.  Since $i$ is an injection,
\[
\alpha_{A} \, | \, D = \alpha_{B} \, | \, D + w
\]
where $w \in \wedge_{GL}(D)$.

\vspace{0.3cm}

\underline{\textbf{Case I}}: \quad $w = d\rho$

Extend $\rho$ to all of $A$, and set $\tilde{\alpha}_{A} = \alpha_{A}
+ d\rho$.  Thus $\{\tilde{\alpha}_{A}\} = \{\alpha_{A}\}$, and
$\tilde{\alpha}_{A} \, | \, D = \alpha_{B} \, | \, D$.  The latter
equation implies there is a unique $\alpha \in \wedge^{odd}(X)$ with
$\alpha \, | \, A = \tilde{\alpha}_{A}$ and $\alpha \, | \, B =
\tilde{\alpha}_{B}$.  Thus by $\ast)$
\begin{eqnarray*}
\bar{\mu} \, | \, A & = & \mu_{A} + i(\{\alpha\}) \, | \, A \\
\bar{\mu} \, | \, B & = & \mu_{B} + i(\{\alpha\}) \, | \, B
\end{eqnarray*}
which implies that $\mu = \bar{\mu} - i(\{\alpha\})$ satisfies the
conditions of the theorem.

\vspace{0.3cm}

\underline{\textbf{Case II}}: \quad $w = g^{\ast}(\Theta) + d\rho$,
where $g : D \rightarrow GL$, and $g^{\ast}(\Theta)$ is not exact.

Using the clutching construction, we may construct a vector bundle $V$
over $X$ with the properties that $V \, | \, A$ and $V \, | \, B$ are
each trivialized by cross-sections $\{E^{A}_{i}\}$ and $\{E^{B}_{i}\}$,
and
\[
\ast\ast) \qquad E^{B}_{j} \, | \, D = \sum_{i} g_{ij} E^{A}_{i} \, |
\, D.
\]
Choose a connection, $\nabla^{\prime}$, on $V$, set $\mathcal{V} = [V,
\{\nabla^{\prime}\}]$ and $\mu^{\prime} = \mathcal{V} - [\textrm{dim}(\mathcal{V})] \in \hat{K}(X)$.

By construction, $\delta(\mu^{\prime} \, | \, A) = 0 =
\delta(\mu^{\prime} \, | \, B)$, and thus
\begin{eqnarray*}
\mu^{\prime} \, | \, A & = & i(\{\alpha^{\prime}_{A}\}) \\
\mu^{\prime} \, | \, B & = & i(\{\alpha^{\prime}_{B}\}) 
\end{eqnarray*}
where $\alpha^{\prime}_{A}, \, \alpha^{\prime}_{B} \in
\wedge^{odd}(A), \, \wedge^{odd}(B)$ \,and\, $\{\alpha^{\prime}_{A}\}$,
$\{\alpha^{\prime}_{B}\}$ represent their equivalence classes modulo $\wedge_{GL}(A)$, $\wedge_{GL}(B)$.

Let $\nabla^{AF}$ and $\nabla^{BF}$ be the flat connections on $V \, |
\, A$ and $V \, | \, B$ defined by making $\{E^{A}_{i}\}$,
$\{E^{B}_{i}\}$ parallel.  By the definition of $i$, and working mod
exact, we may take
\begin{eqnarray*}
\alpha^{\prime}_{A} & = & \mathit{CS}(\nabla^{AF}, \nabla^{\prime} \, | \, A) \\
\alpha^{\prime}_{B} & = & \mathit{CS}(\nabla^{BF}, \nabla^{\prime} \, | \, B).
\end{eqnarray*}
Now, continuing to work mod exact,
\[
\alpha^{\prime}_{A} \, | \, D - \alpha^{\prime}_{B} \, | \, D =
\mathit{CS}(\nabla^{AF} \, | \, D, \nabla^{\prime} \, | \, D) -
\mathit{CS}(\nabla^{BF} \, | \, D, \nabla^{\prime} \, | \, D) =
\mathit{CS}(\nabla^{AF} \, | \, D, \nabla^{BF} \, | \, D) = g^{\ast}(\Theta)
\]
by $\ast\ast)$ and the argument of Lemma 2.3.

Thus, by taking $\bar{\bar{\mu}} = \bar{\mu} - \mu^{\prime}$ and
referring to $\ast)$ we see
\begin{eqnarray*}
\bar{\bar{\mu}} \, | \, A & = & \mu_{A} + i(\{\alpha_{A} -
\alpha^{\prime}_{A}\}) \\
\bar{\bar{\mu}} \, | \, B & = & \mu_{B} + i(\{\alpha_{B} -
\alpha^{\prime}_{B}\})
\end{eqnarray*}
and
\begin{eqnarray*}
(\alpha_{A} - \alpha^{\prime}_{A}) \, | \, D & = & (\alpha_{B} -
\alpha^{\prime}_{B}) \, | \, D + \textrm{exact}.
\end{eqnarray*}
The problem is now reduced to Case I.  \hspace{1 cm}  $\blacksquare$

\vspace{0.3cm}

\underline{\textbf{Corollary 3.10}}: \quad  $\textrm{ker}(ch)$ also
satisfies the Mayer-Vietoris property.

\vspace{0.3cm}

\underline{\textbf{Proof}}: \quad  In the theorem above, if $ch(\mu_{A}) = 0
= ch(\mu_{B})$, then $ch(\mu) \, | \, A = 0 = ch(\mu) \, | \, B$.
Since $ch(\mu)$ is a differential form, this implies $ch(\mu) = 0$.   
 \hspace{1 cm} $\blacksquare$

\vspace{0.3cm}

\underline{\textbf{Proposition 3.11}}: \quad  $\textrm{ker}(ch)$ is a
homotopy functor.

\vspace{0.3cm}

\underline{\textbf{Proof}}: \quad  Any element of $\textrm{ker}(ch)$
is of the form $\mathcal{V} - [\textrm{dim}(\mathcal{V})]$, where
$ch(\mathcal{V}) = \textrm{dim}(\mathcal{V})$.  By 1.9) the pull backs
of $\mathcal{V}$ under two smoothly homotopic $C^{\infty}$ maps would
be isomorphic, and so of course would pull backs of $[\textrm{dim}(\mathcal{V})]$. \qquad $\blacksquare$

\vspace{0.3cm}

\section*{\S4. \, ker$(ch)$ has classifying space the homotopy fibre of
\\
\hspace*{\fill} $BGL \stackrel{ch}{\longrightarrow} \Pi^{\infty}_{n=1} K(C,2n)$ \hspace*{\fill}
}

We begin by introducing a relative group denoted $\pi_{n} (BGL,
\textrm{ker}(ch))$ related to the characterization of homotopy fibers
discussed in Appendix A.  This group will consist of equivalence
classes of stable complex vector bundles with $C$-linear
connections over the n-disk $D^{n}$ so that the Chern-Weil form vanishes on
$\partial D^{n}$.  Two of these $(E, \xi)$ and $(E^{\prime},
\xi^{\prime})$ are equivalent if for some stable bundle isomorphism
over $D^{n}$, 
$E \stackrel{p}{\longrightarrow} E^{\prime}$ the $\mathit{CS}$ form
$\mathit{CS}(\xi, p^{\ast}\xi^{\prime})$ which is closed on $\partial
D^{n}$ is already exact on $\partial D^{n}$.  We can add equivalence
classes by direct sum, and these form a group using Theorem 1.15
applied to pairs, denoted 
$\pi_{n} (BGL, \textrm{ker}(ch))$.

\vspace{0.3cm}

\underline{\textbf{Proposition 4.1}}: \quad  
\[
\pi_{n} (BGL, \textrm{ker}(ch)) = \left\{ 
\begin{array}{cl}
0 & n \, \textrm{odd} \\ 
C & n \, \textrm{even} \\
\end{array}
\right.
\]
with the isomorphism given by the class in the cohomology of the
boundary of the Chern-Simons difference form with the flat connection.

\vspace{0.3cm}

\underline{\textbf{Proof}}: \quad Since $D^{n}$ is contractible, bundles
over $D^{n}$ are trivial, the bundle isomorphism $E \stackrel{p}{\longrightarrow} E^{\prime}$
exists and is unique up to homotopy.  The $\mathit{CS}$ form over
$\partial D^{n} \, \mathit{CS}(\xi, p^{\ast} \xi^{\prime})$ is odd and
therefore exact for $n$ odd.  Thus $\pi_{n} (BGL, \textrm{ker}(ch)) =
0$ for $n$ odd.

For $n$ even, the $\mathit{CS}$ form $\mathit{CS}(\xi, p^{\ast}\xi^{\prime})$
is closed and defines an element in $H^{n-1}(\partial D^{n}, C)
\simeq H^{n}(D, \partial D^{n}; C)$.    $(E, \xi)$ and
$(E^{\prime})$ are equivalent iff this class is zero.  All classes
occur in this by Proposition 2.6 applied to the boundary of the $n$-disk.  This proves Proposition 4.1.
 \qquad $\blacksquare$

\vspace{0.3cm}

\underline{\textbf{Proposition 4.2}}: \quad  The functor $\textrm{ker}(ch)$
is naturally equivalent on pointed compact manifolds with corners to
the based homotopy classes of maps into some classifying space $GL
(C/Z)$.

\vspace{0.3cm}

\underline{\textbf{Proof}}: \quad  $\textrm{ker}(ch)$ of a point is zero.
$\textrm{ker}(ch)$ is a homotopy functor satisfying Mayer-Vietoris by
Corollary 3.10.
$\textrm{ker}(ch)$ also sends finite disjoint unions to finite
products.  It follows from Brown's theorem [10] that $\textrm{ker}(ch)$
on compact manifolds with corners has a classifying space which we denote
$GL (C/Z)$. \\ 
\hspace*{\fill} $\blacksquare$ \qquad

\vspace{0.3cm}

\underline{\textbf{Proposition 4.3}}: \quad $GL (C/Z)$ is homotopy
equivalent to the homotopy fibre of the Chern character map
\[
BGL \stackrel{ch}{\longrightarrow} \Pi^{\infty}_{n=1} K(C,2n).
\]
\vspace{0.3cm}
\underline{\textbf{Proof}}: \quad The map
$\textrm{ker}(ch) \stackrel{\delta}{\longrightarrow} \tilde{K}(Z)$
implies a map $GL (C/Z) \stackrel{\delta}{\longrightarrow} BGL$
where $\tilde{K}(Z)$ is the kernel of the restriction to the base
point which is classified by maps into $BGL$.

The composition $GL (C/Z) \stackrel{\delta}{\longrightarrow} BGL
\stackrel{ch}{\longrightarrow} \Pi^{\infty}_{n=1} K(C,2n)$ is
null homotopic from the definition of $\textrm{ker}(ch)$ as
representing structured bundles with Chern character form identically
zero.  In fact we may consider that we are provided with a preferred
homotopy class of null homotopies for this composition $ch \circ
\delta$.

Using this null homotopy gives a map $(BGL) \cup \, \mathrm{cone} \, GL
(C/Z) \longrightarrow \Pi^{\infty}_{n=1} K(C,2n)$
which then gives a map
\[
\pi_{n} (BGL,GL (C/Z)) \stackrel{CH}{\longrightarrow} \pi_{n}(\Pi^{\infty}_{k=1} K(C,2k)).
\]

In Proposition 4.1 and the paragraph before we have interpreted
$\pi_{n} (BGL,GL (C/Z))$ geometrically as the group $\pi_{n}
(BGL,\textrm{ker}(ch))$ and shown $\mathit{CH}$ is an isomorphism.

By Appendix A, this means $GL (C/Z)
\stackrel{\delta}{\longrightarrow} BGL$ is homotopy equivalent to the
homotopy fibre of $BGL \stackrel{ch}{\longrightarrow} \Pi^{\infty}_{n=1}
K(C,2n).$  \qquad $\blacksquare$

\vspace{0.3cm}

Thus we have 

\underline{\textbf{Theorem 4.4}}:  \quad  $\textrm{ker}(ch)$ is naturally
equivalent to $K^{odd} (C/Z) = $ complex $K$-theory with
coefficients in $C/Z$.

\vspace{0.3cm}

\underline{\textbf{Proof}}: \quad  This will follow from the
definition of $K^{odd}(C/Z)$ and
the above.  It is correct to define $K^{odd} (C/Z)$ as
classified by the homotopy fibre of the map of classifying spaces
corresponding to the map of reduced theories (the kernel of
restrictions to the base points)
\[
\tilde{K}(Z) \stackrel{\otimes C}{\longrightarrow} \tilde{K}(Z)
\otimes_{Z} C \equiv \tilde{K}(C)
\]
i.e.
\[
BGL \stackrel{ch}{\longrightarrow} \Pi^{\infty}_{n=1} K(C,2n)
\]

where $\mathit{ch}$ is the composition of $\otimes C$ with
the Chern equivalence of $\tilde{K}(C)$ and $\Pi^{\infty}_{n=1}
H^{2n}( , C)$.

For then the long exact sequence of the fibration
\[
\textrm{homotopy fibre} \rightarrow BGL \rightarrow \Pi^{\infty}_{n=1} K(C,2n)
\]

becomes
\[
K^{odd}(C) \rightarrow K^{odd}C/Z \rightarrow
K^{even}(Z) \rightarrow K^{even}(C)
\]

or after applying the Chern character equivalence over $C$
\[
H^{odd}(C) \rightarrow K^{odd} C/Z \rightarrow
K^{even}(Z) \rightarrow H^{even}(C).  \qquad \blacksquare
\]

We gather all of this together to arrive at the following result.

\pagebreak

\underline{\textbf{Corollary}}: \quad We have the diagram with exact
diagonals and exact upper and lower boundaries:

\vspace{.5cm}
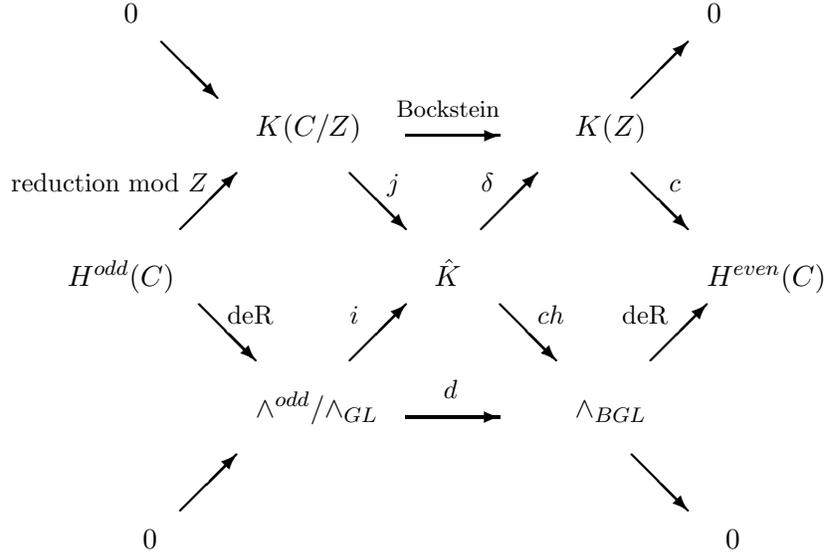
\begin{center}
\setlength{\unitlength}{0.5cm}
\begin{picture}(24,16)\thicklines
\put(5,1){$0$}
\put(20.5,1){$0$}

\put(6,2){\vector(1,1){1.5}}
\put(18,3.5){\vector(1,-1){1.5}}

\put(8,4.5){$\wedge^{odd}/\wedge_{GL}$}
\put(12,4.5){\vector(1,0){2.5}}
\put(16.5,4.5){$\wedge_{BGL}$}
\put(13,5){\small{$d$}}

\put(6.5,7.5){\vector(1,-1){1.5}}
\put(7.25,7){\small{deR}}
\put(10.5,7){\small{$i$}}
\put(10.5,6){\vector(1,1){1.5}}
\put(14.5,7.5){\vector(1,-1){1.5}}
\put(15.5,7){\small{$ch$}}
\put(17.75,7){\small{deR}}
\put(18.5,6){\vector(1,1){1.5}}

\put(3,8){$H^{odd}(C)$}
\put(12.75,8){$\hat{K}$}
\put(20,8){$H^{even}(C)$}

\put(1.5,10.5){\small{reduction mod $Z$}}
\put(6,9.5){\vector(1,1){1.5}}
\put(10.5,11){\vector(1,-1){1.5}}
\put(11.5,10.5){\small{$j$}}
\put(14,10.5){\small{$\delta$}}
\put(14.0,9.5){\vector(1,1){1.5}}
\put(18,11){\vector(1,-1){1.5}}
\put(19,10.5){\small{$c$}}

\put(8,12){$K(C/Z)$}
\put(12,12){\vector(1,0){2.5}}
\put(16.5,12){$K(Z)$}
\put(11.75,12.5){\footnotesize{Bockstein}}

\put(5.5,14.5){\vector(1,-1){1.5}}
\put(18,13){\vector(1,1){1.5}}

\put(4.5,15){$0$}
\put(20,15){$0$}
\end{picture}
\end{center}

\underline{\textbf{Proof}}: \quad The natural equivalence between $\textrm{ker}(ch)$
and $K(C/Z)$ respects the Bockstein sequence because the construction
in Appendix A relates the long exact homotopy sequence of the pair
(total space, fibre) and the long exact sequence of homotopy groups of
a fibration. \hspace{1 cm} $\blacksquare$

\section*{\S5. Hermitian Vector Bundles}

In all that preceded, the the basic objects were complex vector
bundles with connection.  The entire approach immediately applies to
Hermitian bundles with inner product preserving connection.  The same
definition of equivalence goes through and gives rise to a Hermitian
version of Struct.  Analogs of all results remain true, with proofs
following identical lines.

Letting $\hat{K}_{R} = K$ (Hermitian Struct), we obtain the following
commutative diagram,

\vspace{.5cm}
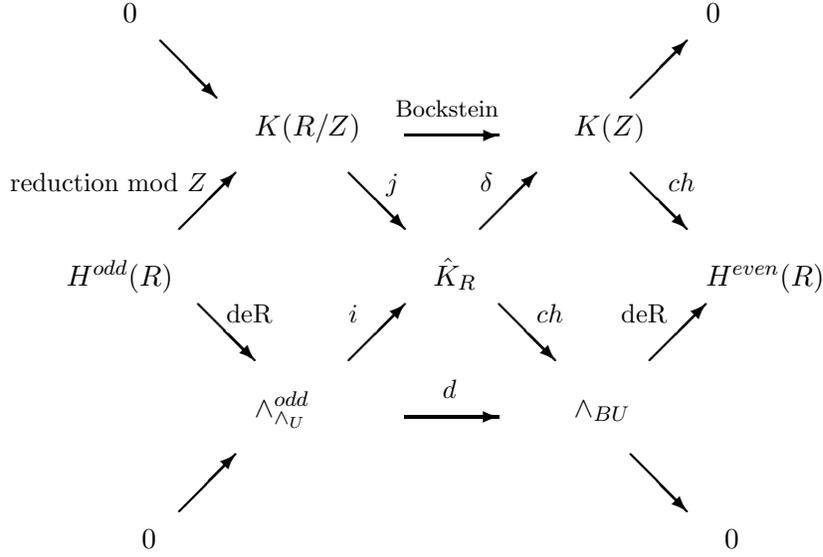
\begin{center}
\setlength{\unitlength}{0.5cm}
\begin{picture}(24,16)\thicklines
\put(5,1){$0$}
\put(20.5,1){$0$}

\put(6,2){\vector(1,1){1.5}}
\put(18,3.5){\vector(1,-1){1.5}}

\put(8,4.5){$\wedge^{odd}_{\wedge_{U}}$}
\put(12,4.5){\vector(1,0){2.5}}
\put(16.5,4.5){$\wedge_{BU}$}
\put(13,5){\small{$d$}}

\put(6.5,7.5){\vector(1,-1){1.5}}
\put(7.25,7){\small{deR}}
\put(10.5,7){\small{$i$}}
\put(10.5,6){\vector(1,1){1.5}}
\put(14.5,7.5){\vector(1,-1){1.5}}
\put(15.5,7){\small{$ch$}}
\put(17.75,7){\small{deR}}
\put(18.5,6){\vector(1,1){1.5}}

\put(3,8){$H^{odd}(R)$}
\put(12.75,8){$\hat{K}_{R}$}
\put(20,8){$H^{even}(R)$}

\put(1.5,10.5){\small{reduction mod $Z$}}
\put(6,9.5){\vector(1,1){1.5}}
\put(10.5,11){\vector(1,-1){1.5}}
\put(11.5,10.5){\small{$j$}}
\put(14,10.5){\small{$\delta$}}
\put(14.0,9.5){\vector(1,1){1.5}}
\put(18,11){\vector(1,-1){1.5}}
\put(19,10.5){\small{$ch$}}

\put(8,12){$K(R/Z)$}
\put(12,12){\vector(1,0){2.5}}
\put(16.5,12){$K(Z)$}
\put(11.75,12.5){\footnotesize{Bockstein}}

\put(5.5,14.5){\vector(1,-1){1.5}}
\put(18,13){\vector(1,1){1.5}}

\put(4.5,15){$0$}
\put(20,15){$0$}
\end{picture}
\end{center}
where $\wedge_{U}$ and $\wedge_{BU}$ are real valued forms, defined
analogously to $\wedge_{GL}$ and $\wedge_{BGL}$.

\vspace{0.3cm}

\underline{\textbf{Corollary 5.1}}:

For any bundle over a closed Riemannian manifold after stabilizing, there is a
unitary connection on the bundle whose Chern-Weil form is the harmonic
representative of the Chern character of the bundle.  Moreover, when
the odd Betti numbers vanish, this structured bundle is unique up to
adding factors with trivial holonomy.

\section*{Appendix A}

Recall in the homotopy theory of spaces homotopy equivalent to $CW$
complexes a map $X \rightarrow
Y$ is homotopy equivalent to the projection map of a Serre fibration.
To see this let us assume $X$ and $Y$ are connected.  First replace $X
\stackrel{p}{\rightarrow} Y$ by $X \stackrel{\tilde{p}}{\rightarrow}
\tilde{Y}$ where $\tilde{p}$ is an inclusion by replacing $Y$ by (the
mapping cylinder of $p$) = $X \times I \cup_{\sim} Y$ where $X
\times 1$ is collapsed by $p$ onto its image in $Y$.

Then replace $X$ by $\tilde{X}$ where $\tilde{X}$ is all the paths in
$\tilde{Y}$ that start in $X$.  Then $\tilde{X}$ maps into $\tilde{Y}$
(continue to call it $\tilde{p}$) with the Serre path lifting property by evaluating a
path at its endpoint in $\tilde{Y}$.  Clearly, $\tilde{Y} \sim Y$,
$\tilde{X} \sim X$, and $\tilde{p} \sim p$.

The fibre $F \rightarrow X$ of $X \stackrel{p}{\rightarrow} Y$ is
defined up to homotopy to be the inclusion into $\tilde{X}$ of the
paths in $\tilde{Y}$ starting in $X$ and ending at a specific point $y
\in Y$ (or $\tilde{Y}$).

\vspace{0.3cm}

\underline{\textbf{Question}}: \quad What properties characterize the
homotopy fibre $F \stackrel{i}{\rightarrow} X$ of a map 
$X \stackrel{p}{\rightarrow} Y$?

\vspace{0.3cm}

\underline{\textbf{Proposition}}: \quad  Suppose we have a map
$F^{\prime} \stackrel{i^{\prime}}{\longrightarrow} X$ and further
suppose the composition 
\mbox{$F^{\prime} \stackrel{i^{\prime}}{\rightarrow} X \stackrel{p}{\rightarrow} Y$} is
provided with a null homotopy so
that the induced map of homotopy sets
\[
\pi_{i}(X,F^{\prime}) \rightarrow \pi_{i}(Y,\textrm{base point})
\]
are bijections $i = 1,2,...$.  Then $F^{\prime}
\stackrel{i^{\prime}}{\rightarrow} X$ is homotopy equivalent to the
homotopy fibre $F \rightarrow X$ of $X \stackrel{p}{\longrightarrow} Y$.

\vspace{0.3cm}

\underline{\textbf{Proof}}: \quad  In this proof we assume $X$ and $Y$ are
connected and $p$ is onto $\pi_{1}$.  Thus $F$ is connected and we
assume $F^{\prime}$ is also connected.  By the path lifting property
of Serre fibrations, the null homotopy of the
composition $F^{\prime} \stackrel{i}{\rightarrow} X
\stackrel{p}{\rightarrow} Y$ defines a canonical homotopy class of
maps $F^{\prime} \rightarrow F$ so that 
\[ \begin{array}{ccccc}
F^{\prime} & \stackrel{i^{\prime}}{\rightarrow} & X & \stackrel{p}{\rightarrow} & Y \\
\downarrow &             & || &             & || \\
F          & \stackrel{i}{\rightarrow} & X  & \stackrel{p}{\rightarrow} & Y \\
\end{array} \]
is homotopy commutative.

Now we look at the exact sequence of homotopy groups and sets
\[ \begin{array}{ccccccccccccccc}
\cdots & \rightarrow & \pi_{2} X & \rightarrow & \pi_{2}(X,F^{\prime}) & \rightarrow &
\pi_{1} F^{\prime} & \rightarrow & \pi_{1} X & \rightarrow &
\pi_{1}(X, F^{\prime}) & \rightarrow & \pi_{0} F^{\prime} & \cong &
\textrm{pt} \\
& & || & & \downarrow & & \downarrow & & || & & \downarrow & & \cong
\downarrow & & \\
\cdots & \rightarrow & \pi_{2} X & \rightarrow & \pi_{2}(X,F) & \rightarrow &
\pi_{1} F & \rightarrow & \pi_{1} X & \rightarrow &
\pi_{1}(X, F) & \rightarrow & \pi_{0} F & \cong &
\textrm{pt} \\
\end{array} \]
In a fibration the Serre path lifting implies the homotopy sets
$\pi_{i}(X,F)$ are isomorphic to $\pi_{i}(Y,\textrm{base point})$ and
thus become groups.  By the above commutative diagram the maps
$\pi_{i}(X,F^{\prime}) \rightarrow \pi_{i}(X,F)$ are bijections.  Thus
the proposition follows from the 5-lemma.   \qquad $\blacksquare$

\section*{References}

\begin{enumerate}
\item Lott, John.  ``$R/Z$ Index Theory''.
Comm Anal Geom 2. 1994. pp. 279-311.

\item Karoubi, M. 
\begin{enumerate}
\item ``Homologie Cyclique et $K$-Th\'{e}orie''.  Asterisque 149. 1987.
\item ``Th\'{e}orie G\'{e}n\'{e}rale des Classes
Caract\'{e}ristiques Secondaires''.  $K$-Theory 4. 1987.  pp. 55-87.
\end{enumerate}

\item Bismut, J.M. 
\begin{enumerate}
\item ``The Index Theorem for Families of Dirac Operators:  Two Heat
Equation Proofs''.  Inv. Math.  83.  pp. 91-151.  1986.
\item with Jeff Cheeger. ``$\eta$-Invariants and Their Adiabatic
Limits''.  J. Amer. Math. Soc. 2.  1989.  pp. 33-70.
\item with D Freed. ``The Analysis of Elliptic Familes II''. 
Comm. Math. Phys. 107.  1986.  pp. 103-163.
\end{enumerate}

\item Bunke, U. and Schick, T. ``Smooth $K$-Theory''.  arXiv
: 0707.0046 (math.KT).  2007. 57 pp.

\item Hopkins, M.J. and Singer, I.M. ``Quadratic functions in
 Geometry, Topology, and $M$-theory''.  J. Diff. Geom 70. 2005. pp. 329-452.

\item Freed, Daniel S.
\begin{enumerate}
\item ``Pions and Generalized Cohomology''.  arXiv: hep-th/0607134.  2006.
29 pp.

\item ``Dirac Charge Quantization and Generalized Differential
Cohomology''.  arXiv:hep-th/0011220. (hep-th,math.DG). 2000.  62 pp.

\item with Gregory W. Moore.  ``Setting the Quantum Integrand of
$M$-Theory''.  \\
arXiv:hep-th/0409135.  2004. 52 pp.  \\
Comm. Math. Physics 263.  2006.  pp. 89-132.

\item with M.J. Hopkins.  ``On Ramond Ramond fields and $K$-Theory''.
arXiv': hep-th/0002027v3. 2000.  12 pp.

\end{enumerate}

\item Simons, James and Sullivan, Dennis.  ``An Axiomatic
Characterization of Ordinary Differential Cohomology''. \\
arXiv : math/0701077v1. (math.AT). 2007. 15 pp. \\
Journal of Topology 1 (1). 2008. pp 45-56.

\item Moore, Gregory W.  Oral Communication.

\item Narasimhan, M.S. and Ramanan, S.  ``Existence of Universal
Connections''.  Amer. J. Math. 83.  1961.  pp. 563-572, 85. 1963. pp. 223-231.

\item Adams, J.F.  ``A Variant of E.H. Brown's Representability
Theorem''.  Topology 10.  1971.  pp. 185-198.  MR 44:1018.

\item Cheeger, Jeff and Simons, James.  ``Differential Characters and Geometric Invariants''.
Notes of Stanford Conference 1973, Lecture Notes in Math. No. 1167.
Springer-Verlag, New York.  1985.  pp. 50-90.

\end{enumerate}

\end{document}